\documentclass[oneside,draft,10pt]{amsart}

\newcommand{\N}{\ensuremath{\mathbb{N}}}
\newcommand{\Z}{\ensuremath{\mathbb{Z}}}
\newcommand{\R}{\ensuremath{\mathbb{R}}}
\newcommand{\dif}{{\rm d}}
\newcommand{\crit}{{\rm Crit}}

\newtheorem{df}{Definition}[section]
\newtheorem{thm}[df]{Teorema}
\newtheorem{rema}[df]{Observaci\'on}
\newtheorem{lema}[df]{Lema}

\usepackage[spanish]{babel}
\usepackage[latin1]{inputenc}
\usepackage{times}
\usepackage{dsfont}
\usepackage[all]{xy}

\title[Homolog\'ia de Morse]{Homolog\'ia de Morse en Variedades Compactas}

\author[C. Marín]{Carlos Alberto Marín Arango}

\address{Departamento de Matem\'atica,\hfill\break\indent  Universidad de Antioquia}
\email{camara@matematicas.udea.edu.co}


\date{May 2010}

\begin{document}

\makeatletter
\renewenvironment{proof}[1][\proofname]{\par
  \pushQED{\qed}%
  \normalfont \topsep6\p@\@plus6\p@\relax
  \trivlist
  \item[\hskip\labelsep
        \scshape
    #1\@addpunct{.}]\ignorespaces
}{%
  \popQED\endtrivlist\@endpefalse
}
\makeatother

\numberwithin{equation}{section}
\swapnumbers

\begin{abstract}
Given a compact Riemannian manifold $(M,g)$ and a Morse function
$f:M\to\R$ whose gradient flow satisfies the Morse--Smale
condition, (i.e. the stable and unstable manifolds of $f$ intersect
transversely) we construct a chain complex called the
Morse--Witten Complex. Our goal on this paper is show that the
homology of the Morse--Witten complex is isomorphic to the
singular homology of $M$.
\end{abstract}
\maketitle

\begin{section}{Introducci\'on}

Los fundamentos de la Teoría de Morse son originalmente introducidos por {\em Marston Morse\/} en el artículo \cite{Morse}. La idea básica es encontrar invariantes topológicos para los puntos críticos de una función diferenciable y estimar técnicas para determinar el valor de tales invariantes. Una referencia clásica para el estudio de ésta es \cite{Milnor}, allí por medio del estudio de los subniveles cerrados
\[f^a=\{x\in M: f(x)\le a\},\] donde $a\in \R$ y $f$ es una función diferenciable definida en una va\-riedad compacta finito dimensional $M$, se obtiene  información sobre el tipo de homotopía de $M$, más específicamente, si $f:M\to \R$ es una aplicación diferenciable definida en una variedad dife\-renciable compacta, los subniveles $f^a$ no cambian cuando el número $a$ varia en un intervalo $[b,c]$ el cual no contiene valores críticos; esto es mostrado deformando los subniveles cerrados por medio de las líneas de flujo del campo gradiente de $f$. Obviamente, la presencia de un valor crítico en el intervalo $[b, c]$ es una obstrucción para tal argumento, ya que en este caso algunas líneas de flujo del campo $\nabla f$ no realizan su trayectoria desde el nivel $c$ hasta el nivel $b$. Si hay un único valor crítico $a\in ]b,c[$, entonces $f^c$ como espacio topológico se obtiene adjuntando al subnivel $f^b$ una célula por cada punto crítico en $f^{-1}(a)$ de dimensión igual al índice del punto crítico. Dado que el pegado de células produce un efecto en la homología de un espacio topológico, la presencia y la cantidad de puntos críticos de un índice específico puede determinarse considerando los grupos de homología de $M$. 

En los últimos 60 años la Teoría de Morse ha estado presente en diversos trabajos, en la decada del 50, Raoul Bott empleando métodos de la teoría de Morse estudió los grupos de homología y homotopía para los espacios simétricos compactos \cite{Bott, Bott1}, de este trabajo se obtiene la prueba del {\em Teorema de periodicidad de Bott}, además se introducen las funciones de Morse--Bott las cuales son una generalización de las funciones de Morse. 
Durante la decada de los 60, la Teoría de Morse es empleada para estudiar algunos aspectos topológicos en variedades, en particular en los trabajos de Stephen Smale los cuales le conducen a la solución de la conjetura de Poncaire para dimensiones mayores a $4$. 
En los años 80, aparece un nuevo enfoque de la teoría de Morse debido a Witten \cite{Witten}, la idea fundamental de este trabajo es asociar a una función de Morse $f:(M,g)\to \R$ definida en una variedad Riemanniana compacta de dimensión finita $n$ un complejo de cadena llamado {\em el complejo de Morse--Witten}, para el cual el $k$-ésimo ($k=0,\dots n$) grupo de cadena es el grupo abeliano libre generado por los puntos críticos de índice $k$ de $f$ y cuyo operador bordo 
realiza un conteo algebraico de las líneas de flujo asociadas con el campo gradiente. La homología del complejo de Morse--Witten es conocida como la {\em Homología de Morse}. A comienzos de la decada de los 90, Andreas Floer apoyado en las ideas de Witten introduce la {\em Homología de Floer} que es la versión infinito dimensional de la Homología de Morse \cite{Floer}. 

Previo al trabajo de Witten, la idea de asociar un complejo a una función de Morse definida en una variedad Riemanniana $M$ fue considerada por René Thom, quien encontró una des\-composición celular para $M$ asociada con $f$; luego Smale introduce una condición adicional en la métrica $g$, de modo que la descomposición celular resultante es un complejo celular.

El objetivo de este trabajo es mostrar que el complejo de Morse--Witten asociado con una función de Morse $f:(M,g)\to \R$ es un complejo de cadena isomorfo al complejo de cadena singular de la variedad $M$. Este resultado aparece probado en \cite{Schwarz}; el enfoque presentado en este trabajo es diferente y se realiza desde el punto de vista de los sistemas dinámicos, a saber, via la intersección de las variedades estable e inestable asociadas con los puntos críticos de la función. La selección de este tópico para la elaboración de este artículo es motivada fundamentalmente en la elegancia y el caracter interdisciplinario de la teoría de Morse, la cual, además de ofrecer una colección considerable de teoremas, introduce conceptos y técnicas que se han tornado herramientas útiles para comprender y solucionar problemas matemáticos en diversas áreas. 

\end{section}

\begin{section}{Notaci\'on y Preliminares}

\begin{subsection}{Puntos cr\'iticos y funciones de Morse}

Sea $M$ una variedad diferenciable de dimensi\'on $n$ y sea $f:M\to \R$ una 
funci\'on diferenciable. Un punto $p \in M$ es llamado un {\em
punto cr\'itico de $f$\/} si la aplicaci\'on lineal inducida $\dif f(p) : T_pM
\to \R$ tiene rango nulo, en este caso el valor real $f(p)$ es llamado 
{\em valor cr\'itico de $f$\/}. Denotamos por $\crit (f)$ el conjunto formado
por todos los puntos cr\'iticos de la funci\'on $f$. Si $a \in \R$, {\em el 
conjunto de los puntos cr\'iticos en el nivel\/} $a$ es definido y denotado
por:
\[\crit _a =  \crit (f) \cap f^{-1}(a).
\] Obviamente, $\crit (f)$ y $\crit _a$ son subconjuntos cerrados de $M$ y el
conjunto de valores regulares de $f$ es dado por 
$\R \setminus f\big(\crit (f)\big)$.

Para cada $p \in \crit (f)$, es posible definir una aplicaci\'on bilineal 
sim\'etrica $\dif^{2}f(p)$ en $T_pM$, llamada {\em la forma Hessiana de $f$ en 
$p$\/}. Un punto cr\'itico $p$ es {\em no degenerado\/} si la forma Hessiana 
de $f$ en $p$ es no degenerada. El {\em \'indice de Morse\/} para el punto 
cr\'itico $p\in M$ se define como siendo el \'indice de la forma Hessiana de $f$ en 
$p$. i.e., la dimensi\'on del mayor subespacio en el cual \'esta es definida 
negativa. Decimos que $f$ es una {\em funci\'on de Morse\/} cuando todos sus 
puntos cr\'iticos son no degenerados.

El Lema de Morse afirma que en coordenadas apropiadas entorno a un punto 
cr\'itico no degenerado, la funci\'on $f$ puede ser descrita por una forma 
cuadr\'atica no degenerada. M\'as espec\'ificamente, si $p$ es un punto cr\'itico no 
degenerado de la funci\'on $f$, entonces existen coordenadas locales 
$(y_1, \cdots , y_n)$ en una vecindad $U$ de $p$ con $y_i(p) = 0$ para cada 
$i$ de modo que la identidad: 
\begin{equation}\label{thm:lemamorse}f = f(p) - (y_1)^2 - \cdots - 
(y_k)^2 + (y_{k+1})^2 + \cdots + (y_n)^2
\end{equation} 
se cumple en $U$, donde $k$ denota el \'indice de Morse de $p$ \cite{Milnor}. 
Como consecuencia de la identidad \eqref{thm:lemamorse} se tiene que el 
\'indice de cualquier punto de m\'aximo local (resp. m\'inimo local) es $n$. 
(resp. $0$) Adem\'as, dado que el origen es el \'unico punto cr\'itico de una 
forma cuadr\'atica no degenerada, se tiene que los puntos cr\'iticos no degenerados de una funci\'on diferenciable $f:M\to \R$ son aislados en 
el conjutno $\crit (f)$. En 
particular, si $f:M\to \R$ es una funci\'on de Morse definida en una variedad 
diferenciable compacta, el conjunto $\crit (f)$ es finito.

Es natural indagar sobre la existencia de funciones de Morse definidas en una 
variedad diferenciable arbitraria $M$. La respuesta a esta pregunta es 
afirmativa en el caso compacto. Adem\'as, en la actualidad es posible mostrar 
que dada una funci\'on diferenciable $g:M\to \R$, existe una 
funci\'on de Morse $f:M\to \R$ arbitrariamente cerca de $g$ (respecto
a la topolog\'ia de convergencia uniforme) \cite{Dubrovin2}.
\end{subsection}
\begin{subsection}{Orientaci\'on, Transversalidad y N\'umero de intersecci\'on}
Una {\em orientaci\'on\/} (en el sentido diferenciable) para una variedad diferenciable $M$ es
definida como una aplicaci\'on que asocia a cada punto $x\in M$ una 
orientaci\'on para el espacio tangente $T_xM$; la cual depende continuamente de 
$x$, tal continuidad puede establecerse en t\'erminos de la existencia de un 
atlas formado por cartas positivamente orientadas. En el caso de las 
variedades topol\'ogicas no se posee en general un espacio tangente en cada 
punto, por lo cual, la definici\'on de orientaci\'on presentada para el caso de
una variedad diferenciable no tiene sentido. Sin embargo, empleando teor\'ia 
de homolog\'ia es posible obtener una definici\'on elegante para el concepto 
de orientaci\'on en una variedad, la cual incluye las variedades 
topol\'ogicas. En efecto, una {\em orientaci\'on homol\'ogica\/} para una
variedad topol\'ogica $n$-dimensional $M$ es una secci\'on global del fibrado 
$\mathcal{O}(M)$ sobre $M$, en que
\[
\mathcal{O}(M)=\bigcup_{x\in M}H_n(M,M\setminus \{x\}),
\]
la cual asocia a cada punto $x\in M$ un generador del grupo c\'iclico infinito 
$H_n(M,M\setminus \{x\})$. Para el caso de una variedad 
diferenciable $M$, dado $x\in M$ es posible 
mostra que existe una biyecci\'on can\'onica entre el conjunto de las 
orientaciones en el sentido homol\'ogico de $M$ en el punto $x$ (i.e., el 
conjunto de los generadores del grupo c\'iclico infinito 
$H_n\big(M,M\setminus \{x\}\big)$) y el conjunto de orientaciones en el sentido 
diferenciables en $x$ (i.e., el conjunto de orientaciones para el espacio 
vectorial $T_xM$).  Una {\em orientación transversal\/} para una subvariedad 
$N \subset M$ es una aplicaci\'on que asocia a cada punto $x\in N$ una 
orientaci\'on para el espacio vectorial 
$T_xM/T_xN$, la cual dependa continuamente de $x$. El siguiente resultado es adaptado de ideas encontradas en \cite{Spanier}.
\begin{thm}\label{moustro}
Si $N\subset M$ es una subvariedad cerrada tranversalmente orientada y 
codimensi\'on $n$, existe un isomorfismo:
\[H_n(M,M\setminus N) \xrightarrow{\;\;\varrho\;\;} H_0(N).\] Adem\'as, si
 $M'$ es una subvariedad de $M$ transversal a $N$, y
$N'=M'\cap N$, entonces el siguiente diagrama: 
\begin{equation} \xymatrix{
H_n(M,M\setminus N) \ar[r]^-\varrho& H_0(N)\\
H_n(M',M'\setminus N') \ar[r]^-\varrho
\ar[u]^{i_*} & H_0(N') \ar[u]_{i_*} } 
\end{equation} es conmutativo. En particular, si 
$N=\{a\}$ con ${\rm{dim}} M = n$, entonces el isomorfismo
$$\varrho : H_n(M,M\setminus \{a\}) \longrightarrow H_0(\{a\})$$
mapea el generador asociado con la  orientaci\'on diferenciable en $T_aM$, 
sobre el generador can\'onico de $H_0(\{a\})$.
\end{thm}

Al igual que para el concepto de orientación en variedades, la noción de {\em número de intersección\/} de una función con una suvariedad 
de su codominio puede ser formulada tanto en el lenguaje de la topología diferencial (la cual exige la diferenciabilidad de la función en cuestión) \cite{Pollack}, como en el lenguaje topología algebraica (la cual sólo requiere la continuidad de la función). En el lenguaje de la topología algebraica el caso más simple corresponde a 
una función continua definida en un subconjunto abierto de la esfera $S^n$ a valores en una variedad topológica $n$-dimensional orientada $N$. Denotamos por $\tau^{[n]}$ el generador del grupo $H_n\big(\R^n,\R^n\setminus\{0\}\big)$ el cual es asociado a la orientación canónica de $\R^n$ y por $\alpha^{[n]}$\footnote{la orientación diferenciable asociada a la orientación $\alpha^{[n]}$ of $S^n$ es precisamente la orientación dada por el vector normal que apunta para fuera. \cite{Tausk}.} 
el generador del grupo $\tilde{H}_n(S^n)$ obtenido de $\tau^{[n]}$ 
via el isomorfismo $\partial_*$ que aparece en la secuencia larga en homología del par $\big(\R^n,\R^n\setminus\{0\}\big)$. Sean $U\subset S^n$ un conjunto abierto, $P\subset N$ una subvariedad cerrada transversalmente orientada y $f : U \to N$ una función continua tal que  
$f^{-1}(P)$ es un conjunto compacto; el {\em número de intersección \/} (en el sentidio homológico) de la función  
$f$ con la subvariedad $P$ se define como el entero $\eta(f,P)$ para el cual se cumple la igualdad:
\[
\phi\big(\alpha^{[n]}\big)=\eta(f,P).
\]donde $\phi:\tilde{H}_n(S^n)\to \Z$ es el homomorfismo dado por la composición que aparece en el diagrama:
\[\xymatrix{
\tilde{H}_n(S^n) \ar[r]^-{i_*} \ar[dddrr]_{\phi} & H_n\big(S^n,S^n\setminus f^{-1}(N)\big) & H_n(U,U\setminus f^{-1}(N)) \ar[l]^-{exc}_-{\simeq} \ar[d]^{\ \ f_*}\\
& &  H_n(M,M\setminus N)\ar[d]^{\varrho}\\
             & &  H_0(N) \ar[d]^{\oplus} \\
             & &  \Z } 
\]

\end{subsection}
\end{section}

\begin{section}{El complejo de Morse--Witten}
 
Sean $(M,g)$ una variedad Riemanniana compacta $n$-dimensional y $f:M\to \R$ una funci\'on de Morse. El objetivo
es asociar un complejo de cadena a $f$, o m\'as precisamente a $-\nabla f$, 
en que $\nabla f$ denota el campo gradiente de $f$ respecto a $g$. Claramente la singularidades de $-\nabla f$ son precisamente los puntos
cr\'iticos de $f$; además, por compacidad este es un campo vectorial completo cuyo flujo induce una acci\'on 
diferenciable:
\[
(t,x) \longmapsto t\cdot x
\] del grupo aditivo $\R$ sobre $M$ y para el cual cada l\'inea de flujo converge a un punto cr\'itico 
de $f$, i.e., dado un punto arbitrario $x \in M$ se tiene que 
$\lim_{t \rightarrow \pm\infty} t\cdot x$ existe y es un punto cr\'itico de $f$.
Para este sistema din\'amico, las variedades {\em estable\/} e 
{\em inestable\/} asociadas al punto cr\'itico $p \in M$ se definen 
respectivamente por:
\[\mathcal{W}_s(p)=\{ x \in M : \lim_{t\to +\infty} t\cdot x=p \}
\]
\[ \mathcal{W}_u(p)=\{ x \in M : \lim_{t\to -\infty} t\cdot x=p
\}.\]

\noindent Tanto $\mathcal{W}_s(p)$ como $\mathcal{W}_u(p)$ son
subvariedades de $M$, homeomorfas a los espacios euclídeos $\R^{n-\rm{ind}\_(p)},
\R^{\rm{ind}\_(p)}$ respectívamente; donde $\rm{ind}\_(p)$ denota el \'indice de 
Morse del punto crítico $p$, \cite{Palio}.
Obviamente si $x$ pertenece a la variedad estable (resp., inestable) de un 
punto critico $p$ entonces $t\cdot x$ tambi\'en pertenece a la variedad estable 
 (resp., inestable) de $p$. En particular, para $x\in\mathcal{W}_s(p)$ se tiene:
\begin{equation}\label{eq:gradtang}
-\nabla f(x)=\frac{\dif}{\dif t}\;t\cdot x\Big\vert_{t=0}\in T_x\mathcal{W}_s(p),
\end{equation}
para todo $x\in\mathcal{W}_s(p)$, similarmente $-\nabla f(x)\in T_x\mathcal{W}_u(p)$, para todo $x\in\mathcal{W}_u(p)$. 

Dados puntos $p,q\in \crit (f)$ tales que las variedades $\mathcal{W}_u(p)$,
$\mathcal{W}_s(q)$ son transversales y no dijuntas, como consecuencia 
de la transversalidad se tiene que $\mathcal{W}_u(p)\cap \mathcal{W}_s(q)$ 
es una subvariedad embebida de $M$, cuya dimensi\'on es calculada como sigue:
\[\begin{array}{rcl}
{\rm dim}\big(\mathcal{W}_u(p)\cap\mathcal{W}_s(q)\big)&=&{\rm dim}\big(\mathcal{W}_u(p)\big)+{\rm dim}\big(\mathcal{W}_s(q)\big)-{\rm dim}(M)\\[0.2cm]
&=& \rm{ind}\_(p)+{\rm dim}(M)-\rm{ind}\_(q)-{\rm dim}(M)\\[0.2cm]
&=& {\rm{ind}}\_(p)-{\rm{ind}}\_(q).
\end{array}\]
Adem\'as, cuando $p\ne q$, para cada 
$x\in\mathcal{W}_u(p)\cap\mathcal{W}_s(q)$ se tiene:
\[
T_x\big(\mathcal{W}_u(p)\cap\mathcal{W}_s(q)\big)=T_x\mathcal{W}_u(p)\cap T_x\mathcal{W}_s(q),
\] en particular, $\rm{ind}\_(p)>\rm{ind}\_(q)$ y $-\nabla f(x)\in T_x\big(\mathcal{W}_u(p)\cap\mathcal{W}_s(q)\big)$.

El concepto de transversalidad entre las variedades estable e inestable es introducido en teor\'ia por 
Smale, quien descubri\'o que este requisito para la intersecci\'on de estas va\-riedades se mantiene v\'alido para una m\'etrica 
generica en $M$. Por esta raz\'on, tal condici\'on es conocida como
la condici\'on de Morse-Smale; m\'as espec\'ificamente:
\begin{df}
Dado un n\'umero entero $k$, la función $f:(M,g)\to \R$ satisface la  
{\em condici\'on de Morse--Smale de orden $k$\/} si para cada par de puntos 
cr\'iticos $p,q\in M$ con $\rm{ind}\_(p)-\rm{ind}\_(q)\le k$, la variedad 
inestables de $p$ y la variedad estable de $q$ son transversales.
Si $f:(M,g)\to\R$ satisface la condici\'on de Morse--Smale para cada $k\in\Z$ 
simplemente diremos que $f$ satisface la {\em condici\'on de Morse--Smale\/}. 
\end{df}

La condición de transversalidad de Morse--Smale, permite definir una relación de orden en los índices de los puntos del conjunto $\crit (f)$ determinada por la existencia de una línea de flujo entre dos puntos críticos, de forma que el 
\'indice de Morse decrese a lo largo de las l\'ineas de flujo de $-\nabla f$. 
Adem\'as si $\rm{ind}\_(p)-\rm{ind}\_(q)=1$ y las variedades $\mathcal{W}_u(p), 
\mathcal{W}_s(q)$ son no disjuntas, entonces la variedad
$\mathcal{W}_u(p)\cap\mathcal{W}_s(q)$ tiene dimensi\'on $1$ y consiste de 
una colecci\'on de l\'ineas de flujo para $-\nabla f$ uniendo el punto 
cr\'itico $p$ con el punto cr\'itico $q$. Es posible mostrar que dicha 
colecci\'on de l\'ineas de flujo es finita. A saber, si asumimos que la 
funci\'on de Morse $f:(M,g)\to\R$ satisface la condici\'on de Morse--Smale de 
orden $1$, empleando la compacidad se puede obtener la clausura del conjunto 
$\mathcal{W}_s(p)$ adicionando algunas {\em l\'ineas de 
flujo por pasos\/} para el campo $-\nabla f$, i.e., l\'ineas de flujo que 
pasan por otros puntos cr\'iticos antes de alcanzar su punto final. 
M\'as espec\'ificamente, en \cite{Tausk} es probado el siguiente resultado: 
Sea $p\in\crit (f)$. Si $x\in\overline{\mathcal{W}_s(p)}$ entonces existe una 
$k$-l\'inea de flujo por pasos desde el punto $-\infty\cdot x$ hasta el punto $p$. 
Donde por una $k$-l\'inea por pasos entendemos una sucesi\'on 
$\gamma = (\gamma_1,...,\gamma_k)$ de l\'ineas de flujo $\gamma _i : \R \to M$ 
para $-\nabla f$, tales que $\lim_{t \to +\infty} \gamma _i (t) = \lim_{t \to -\infty} \gamma _{i + 1} (t)$ para cada $i = 1,...,k-1$. 

\begin{lema}
Dados $p,q\in \crit (f)$ con $\rm{ind}\_(p)-\rm{ind}\_(q)=1$, el n\'umero de
 l\'ineas de flujo desde $p$ hasta $q$ es finito.
\end{lema}
\begin{proof}
Sea $a<f(p)$ tal que no hay valores cr\'iticos de $f$ en el intervalo 
$\left[a,f(p)\right[$, cada l\'inea de flujo no constante de $-\nabla f$ 
contenida en $\mathcal{W}_u(p)$ intersecta el nivel $f^{-1}(a)$ (precisamente 
en un punto), o sea, existe una bijecci\'on entre el conjunto de las l\'ineas 
de flujo de $-\nabla f$ desde $p$ hasta $q$ y la subvariedad cero dimensio\-nal
$\mathcal{W}_u(p)\cap\mathcal{W}_s(q)\cap f^{-1}(a)$. Veamos que 
$\mathcal{W}_u(p)\cap\mathcal{W}_s(q)\cap f^{-1}(a)$ es un conjunto finito.  
$\mathcal{W}_u(p)\cap\mathcal{W}_s(q)\cap f^{-1}(a)$ es un conjunto discreto, 
para obtener el resultado deseado es necesario mostrar que es un conjunto 
cerrado. Para esto, sea $x$ un punto en su clausura, como 
$\mathcal{W}_u(p)\cap f^{-1}(a)$ es compacto (homeomorfo a la esfera), se tiene
que $x \in \overline{\mathcal{W}_s(q)}\cap \mathcal{W}_u(p)\cap f^{-1}(a)$. 
Por lo tanto, existe una $k$-l\'inea por pasos desde $q_1 = \lim_{t \to -\infty} t\cdot x \in \crit(f)$ hasta $q$, luego existe una $k+1$-l\'inea por pasos 
desde $p$ hasta $q$. Como $1={\rm{ind}}\_(p)-{\rm{ind}}\_(q)\ge k+1$, 
necesariamente $k=0$ lo que implica $q_1=q$.
\end{proof}

Ahora estamos en condiciones de presentar el {\em Complejo de Morse--Witten\/} 
asociado al campo $-\nabla f$. Antes de proceder, para cada par de puntos
$p, q \in \crit (f)$ es necesario escoger de forma apropiada una 
orientaci\'on en $\mathcal{W}_u(p)\cap\mathcal{W}_s(q)$. En efecto, dado $p\in \crit (f)$ escojemos una 
orientaci\'on para el espacio vectorial $T_p\mathcal{W}_u(p),$ esta selecci\'on
induce una orientaci\'on en la variedad inestable $\mathcal{W}_u(p)$; 
adem\'as una orientaci\'on transversal en la va\-riedad es\-table 
$\mathcal{W}_s(p)$. Por lo tanto, dados $p,q\in \crit (f)$ con 
$\rm{ind}\_(p)-\rm{ind}\_(q)=1$. En cualquier 
punto sobre una l\'inea de flujo desde $p$ hasta $q$ se tiene un isomorfismo
can\'onico: 
\begin{equation}\label{eq:exactayorientacion} 
T\mathcal{W}_u(p) \to T\big(\mathcal{W}_u(p)\cap \mathcal{W}_s(q)\big) \oplus \big(TM/T\mathcal{W}_s(q)\big) 
\end{equation}
La orientaci\'on en $\mathcal{W}_u(p)\cap \mathcal{W}_s(q)$ es escogida de 
forma tal que el isomorfismo \eqref{eq:exactayorientacion} sea positivo.

\begin{df}{(El complejo de Morse--Witten)}
Para $k\ge 0$, se define el grupo de cadena $\mathcal{C}_k(f)$ como siendo el 
grupo libre abeliano generado por el conjunto de puntos cr\'iticos de $f$ con 
\'indice de Morse $k$; y el operado bordo 
$\partial _k:\mathcal{C}_k(f) \to \mathcal{C}_{k-1}(f)$ se define por: 
\begin{equation}\label{eq:bordo}
 \partial_{k}(p) =\sum_{q\in \crit_{k-1}(f)}
 n(p,q)q.
\end{equation}
Donde el coeficiente para $q\in \crit_{k-1}(f)$ en la expresi\'on para el bordo
de $p$, es dado por un conteo algebraico del n\'umero de l\'ineas de flujo 
contenidas en $\mathcal{W}_u(p)\cap\mathcal{W}_s(q)$. El signo de cada 
l\'inea de flujo es determinado comparando la orientaci\'on natural inducida
por el campo vectorial $-\nabla f$ con la orientaci\'on inducida por 
\eqref{eq:exactayorientacion}. 
La {\em Homolog\'ia de Morse\/} para $M$ se define como siendo la homolog\'ia 
del complejo $(\mathcal{C}(f)_*; \partial)$. 
\end{df}
 
El proposito de este trabajo es demosrar que el complejo de Morse--Witten es un complejo de cadena cuyos grupos de homología son isomorfos a los grupos de homología singular de $M$.

\begin{rema}
Dados $p,q \in \crit (f)$ con $n\_(p) - n\_(q) = 1$, $a \in \R$ de modo que 
 no hay valores cr\'iticos para $f$ en el intervalo $\big(f(q),f(p)\big)$, 
el coeficiente para $q$ en la expresi\'on para el bordo del punto $p$, 
 que aparece en \eqref{eq:bordo}, coincide con el n\'umero
de intersecci\'on de las esferas $f^{-1}(a) \cap \mathcal{W}_u(p)$, $f^{-1}(a) \cap \mathcal{W}_s(q)$ en $f^{-1}(a)$. En efecto, el n\'{u}mero de
intersecci\'on de tales esferas es dado por el conteo algebraico de los puntos en el conjunto finito
$\mathcal{W}_u(p) \cap \mathcal{W}_s(q) \cap f^{-1}(a)$, en que 
el signo de un punto $x$ es positivo, si el isomorfismo:
\[
T_x\big(\mathcal{W}_u(p) \cap f^{-1}(a)\big) \xrightarrow{\;\pi\dif
i(x)\;} \frac{T_{x}\big(f^{-1}(a)\big)}{T_{x}\big(\mathcal{W}_s(q)
\cap f^{-1}(a)\big)}
\] es positivo; caso contrario, es negativo. Por otro lado, la aplicaci\'on $T_x\big(f^{-1}(a)\big)
\xrightarrow{\;\;\pi\;\;} T_xM/T_x\big(\mathcal{W}_s(q)\big)$ induce un 
isomorfismo:
\begin{equation}\label{estrella}
\frac{T_x\big(f^{-1}(a)\big)}{T_x\big(f^{-1}(a) \cap
\mathcal{W}_s(q)\big)} \xrightarrow{\;\;\simeq \;\;}
\frac{T_xM}{T_x\big(\mathcal{W}_s(q)\big)}. 
\end{equation} Como el espacio vectorial $T_xM/T_x\big(\mathcal{W}_s(q)\big)$ 
es orientado, empleando el isomorfismo (\ref{estrella}) es posible determinar 
cuando una base para $T_x(f^{-1}(a))/T_x\big(\mathcal{W}_s(q)\big)$ es positiva o
negativa. A saber, si $(b_1, \dots ,
b_{k-1})$ es una base para el espacio $T_x\big(\mathcal{W}_u(p) \cap f^{-1}(a)\big)$
tal que $(-\nabla f(x),b_1, \dots , b_{k-1})$ sea una base positiva
de $T_x\big(\mathcal{W}_u(p)\big)$. Luego la l\'inea de flujo $\sigma$
pasando por el punto $x \in \mathcal{W}_u(p) \cap
\mathcal{W}_s(q)$, tiene signo positivo si, y solamente si, $[b_1],
\dots , [b_{k-1}]$ es una base positiva de
$T_xM/T_x\big(\mathcal{W}_s(q)\big)$. Por otro lado, si $(b_1,
\dots , b_{k-1})$ es una base positiva para $T_x\big(\mathcal{W}_u(p)
\cap f^{-1}(a)\big),$ entonces un punto $x \in \mathcal{W}_u(p) \cap
\mathcal{W}_s(q)$ tiene signo positivo si, y solamente si, 
$[b_1], \dots , [b_{k-1}]$ es una base positiva para el espacio
$T_x(f^{-1}(a))/T_x\big(\mathcal{W}_s(q)\cap f^{-1}(a)\big)$ si, y solamente si, $[b_1], \dots , [b_{k-1}]$ es una base positiva para el espacio
$T_xM/T_x\big(\mathcal{W}_s(q)\big)$.
\end{rema}


Las funciones $f:M\to \R$ 
cuyos valores cr\'iticos est\'an ordenados de acuerdo a los \'indices de los 
correspondientes puntos cr\'iticos, i.e., $f(p)=f(q)$ siempre que los 
respectivos \'indices ${\rm{ind}}\_(p)$, ${\rm{ind}}\_(q)$ sean iguales, y
$f(p)>f(q)$ si ${\rm{ind}}\_(p)>{\rm{ind}}\_(q)$ son llamadas {\em Funciones de Smale.\/} A diferencia de las
funciones de Morse, el conjunto de las funciones de Smale no necesariamente es 
un conjunto denso en el espacio de las funciones continuas $M\to \R$, 
\cite{Dubrovin3}. Sin embargo, se tiene el siguiente resultado adaptado de 
ideas encontradas en \cite{Milno}.

\begin{thm}\label{thm:autoindex}
Sea $f : M \to \R$ una funci\'on de Morse satisfaciendo la condici\'on de 
Morse-Smale de orden cero. Entonces existen una funci\'on de Morse $\tilde{f} : M \to \R$ y una m\'etrica Riemanniana 
$\tilde{g}$ en $M$ tales que el gradiente de $f$ respecto a 
$g$ coincide con el gradiente de $\tilde{f}$ respecto a 
$\tilde{g}$. 
\end{thm}

Actualmente es posible mostrar que la función $\tilde{f}$ del teorema anterior es {\em auto--indexante\/}, i.e., satisface la condición adicional $\tilde{f}(p)={\rm{ind}}\_(p)$ para cada punto crítico $p$.
Dado que el complejo de Morse--Witten asociado a la función $f:(M,g)\to \R$ depende exclusivamente del campo  
gradiente $\nabla f$, con la notaci\'on del teorema anterior, el complejo de Morse--Witten asociado con
la funci\'on $f:(M,g)\to \R$ coincide con el complejo de Morse--Witten asociado a la funci\'on 
$\tilde{f}:(M,\tilde{g})\to \R$. Por lo tanto, sin p\'erdida de generalidad con el fin de estudiar la homología de Morse en la variedad $M$, 
es posible asumir que $f : M \to \R$ es una funci\'on de Morse Auto--indexante
que satisface la condici\'on de Morse-Smale de orden $1$.

Sea $k>0$, bajo la suposiciones anterioriores, el subnivel $f^{k-1}$ es un 
retrato por deformaci\'on (fuerte) del conjunto $f^k\setminus \crit_k$. Por lo tanto, se tiene un isomorfismo en homolog\'ia
\begin{equation}\label{eq:excisioniso1}
H\big(f^k,f^{k-1}\big)\cong H\big(f^k,f^{k}\setminus \crit_k\big).
\end{equation} 
inducido por inclusi\'on. Si $\crit_k=\{p_1,p_2, \dots, p_s\}$ y escogemos 
conjuntos abiertos disjuntos $(U_i)_{i=1}^s$ en $M$ tales que $p_i\in U_i$, 
$i=1,\dots, r$. Sean $U=\cup_{i=1}^sU_i$ y $U^*=\cup_{i=1}^sU_i\setminus\{p_i\}$. 
Como $\crit_k$ es un subconjunto cerrado del abierto relativo $U\cap f^k$ de 
$f^k$, como consecuencia del principio de escisi\'on, obtenemos el siguiente 
isomorfismo en homolog\'ia inducido por inclusi\'on:
\begin{equation}\label{eq:excisioniso2}
H\big(f^k\cap U,f^k \cap U^*\big)\cong H\big(f^k,f^k\setminus\crit_k\big).
\end{equation} Adem\'as,
\begin{equation}\label{eq:excisioniso3}
\begin{array}{ccl}
H\big(f^k\cap U,f^k\cap U^*\big) & \cong & \displaystyle \bigoplus_{i=1}^{s}
H\big(f^{k}\cap U_i,(f^{k}\cap U_i)\setminus \{p_i\}\big)\\ [0.5cm] &
\cong & \displaystyle \bigoplus_{i=1}^{s} H\big(f^{k},f^{k}\setminus \{p_i\}\big).
\end{array}
\end{equation}
Como consecuencia de los isomorfismos \eqref{eq:excisioniso1}, 
\eqref{eq:excisioniso2}, \eqref{eq:excisioniso3} se obtiene un isomorfismo
\begin{equation}\label{eq:excisioniso4}
H\big(f^{k},f^{k-1}\big)\cong \bigoplus_{i=1}^{s} H\big(f^{k},f^{k}\setminus
\{p_i\}\big). 
\end{equation} 
Por el lema de Morse, existe un conjunto abierto entorno de cada punto 
$p_i\in \crit_k$ en $f^k$ el cual es homeomorfo a una vecindad del origen en 
el cono:
\[C=\{(x,y)\in \R^k\times \R^{n-k}: \left\|x\right \|\ge \left\|y\right\|\},\]
por medio de un homeomorfismo que mapea el punto $p_i$ sobre el origen. Luego
\begin{equation}\label{eq:excisioniso5}
H(f^{k},f^{k}\setminus \{p_i\})\cong H(C,C\setminus \{0\})\cong H(\R^k,\R^k\setminus\{0\}),
\end{equation} Por lo tanto, para cada $p\in \crit_k$ 
\[H_i(f^{k},f^{k}\setminus \{p\})\cong \left\{ \begin{array}{ll}
\displaystyle \Z,  &  i = k \\
\displaystyle 0, & i\ne k
\end{array} \right.
\]Sigue de \eqref{eq:excisioniso4} que
\[ H_i(f^{k},f^{k-1}) \cong \left\{
\begin{array}{ll}
\displaystyle \oplus_{\lvert \crit _k \rvert}\Z, & i=k \\
\displaystyle 0, & i\ne k
\end{array} \right.
\]Como consecuencia se tiene que los subniveles cerrados 
$(f^k)_{k\ge0}$ asociados a la función $f$ forman un buena filtraci\'on para la variedad $M$. Por lo 
tanto, la homolog\'ia del complejo de cadena asociado con esta filtraci\'on es 
isomorfa con la homolog\'ia singular de $M$.

\begin{rema}
Como consecuencia del isomorfismo \eqref{eq:excisioniso5}, se tiene
\begin{equation}\label{eq:excisioniso6}
H(\mathcal{W}_u(p_i),\mathcal{W}_u(p_i)\setminus \{p_i\})\cong H(\R^k,\R^k\setminus\{0\})\cong  H(f^k,f^k\setminus \{p_i\}).
\end{equation}
\end{rema}
\end{section}

\begin{section}{Teorema Fundamental}

Dedicamos esta secci\'on a mostrar el resultado principal de este trabajo, 
a saber:
\begin{thm}\label{thm:principal}
Sea $f:M\to \R$ una funci\'on de Morse auto--indexante, definida en una variedad
Riemanniana compacta $(M,g)$ satisfaciendo la condici\'on de Morse--Smale de 
orden $1$. Entonces el complejo de Morse--Witten $(\mathcal{C}_k(f),\partial _k)_{k \in Z}$ es un complejo de cadena, el cual es isomorfo al complejo de cadena singular de $M$ \cite{Munkres}.
\end{thm}

Como consecuencia de lo observado en la sección anterior, para obtener el 
resultado deseado, es suficiente probar que el complejo de cadena asociado 
con la filtraci\'on $(f^k)_{k\ge 0}$ de $f$, es isomorfo al complejo de 
Morse--Witten de $f$. Para esto, consideremos fija una funci\'on $f:M\to \R$ 
como en el enunciado del teorema \ref{thm:principal}. Para cada $k\ge 0$ 
construiremos un isomorfismo: $\rho:\mathcal C_k(f)\to H_k(f^k,f^{k-1})$ tal 
que el diagrama:
\begin{equation}\label{d0}
\xymatrix{\mathcal{C}_k(f) \ar[r]^-{\rho} \ar[d]_{\partial_k} & H_k(f^k,f^{k-1}) \ar[d]^{\partial _*}\\
\mathcal{C}_{k-1}(f) \ar[r]_-{\rho}  & H_{k-1}(f^k,f^{k-1}) } 
\end{equation}
sea conmutativo. 

\noindent Sea $k\ge 0$ un número entero, empleando los isomorfismos \eqref{eq:excisioniso4} y \eqref{eq:excisioniso6}
 se define un isomorfismo:
\[\bigoplus_{p_i\in \crit_k} H_k\left(\mathcal{W}_u(p_i),\mathcal{W}_u(p_i)^*\right)\cong H_k(f^k,f^{k-1})\]
requiriendo la conmutatividad del diagrama:
\begin{equation}\label{quaserho}
\xymatrix{%
\displaystyle\bigoplus_{p_i\in\crit_k}H_k\big(\mathcal W_u(p_i),\mathcal
W_u(p_i)^*\big)\ar[rd]_-{\varphi}
\ar[r]^-{\scriptscriptstyle\cong\;}&
\displaystyle\bigoplus_{p_i\in\crit_k}H_k\big(f^k,f^k\setminus\{p_i\}\big)\\
&H_k(f^k,f^{k-1})\ar[u]_-{\scriptscriptstyle \cong\;\;\eqref{eq:excisioniso4}}}
\end{equation} donde empleamos la notaci\'on 
$\mathcal{W}_u(p)^*=\mathcal{W}_u(p)\setminus \{p\}$.

Para cada $p\in \crit (f)$, escogiendo una orientaci\'on para el espacio 
vectorial $T_p\mathcal W_u(p)$, se obtiene un generador $\mathcal O_p$ 
para el grupo c\'iclico infinito
$H_k\big(\mathcal W_u(p),\mathcal{W}_u(p)^*\big)$; y por lo tanto un 
isomorfismo:
\begin{equation}\label{isoCkHk} 
\mathcal
C_k(f)\xrightarrow{\;\;\simeq\;\;}\bigoplus_{p_i\in\crit_k}H_k\big(\mathcal
W_u(p_i),\mathcal W_u(p_i)^*\big). 
\end{equation} Definimos $\rho$ 
como siendo el isomorfismo obtenido por la composici\'on del isomorfismo 
\eqref{isoCkHk} con el isomorfismo $\varphi$ que aparece en el diagrama 
\eqref{quaserho}. 
Con la notaci\'on y terminolog\'ia anterior, la conmutatividad del 
diagrama \eqref{d0} se obtiene si para cada 
$p\in \crit_k$, cada $q\in \crit_{k-1}$ vale la igualdad:
\[
\Hat{q}\big(\partial_k(p)\big)=\Hat{q}\big(\rho ^{-1}\circ \partial_*\circ \rho (p)\big)
\] donde $\Hat{q}:\mathcal
C_{k-1}(f)\to \Z$ es la aplicaci\'on que da el coeficiente 
en $q$.

\noindent Para cada $q \in \crit _{k-1}$, como consecuencia de la definici\'on 
de $\rho$ y de la conmutatividad del diagrama: 
\[
\xymatrix{\displaystyle\bigoplus_{q_i\in \crit_{k-1}}H_{k-1}\big((\mathcal W_u(q_i),\mathcal{W}_u(q_i)^*\big) \ar[r]^-{\cong} & \displaystyle\bigoplus_{q_i\in \crit_{k-1}}H_{k-1}(f^{k-1},f^{k-1}\setminus \{q_i\})\\
H_{k-1}\big((\mathcal W_u(q),\mathcal{W}_u(q)^*\big) \ar[u]^-{i_*} \ar[r]_{\cong} & H_{k-1}(f^{k-1},f^{k-1}\setminus \{q\})\ar[u]_-{i_*} }  
\]se tiene que el diagrama: 
\begin{equation}\label{d2}
\xymatrix{H_{k-1}(f^{k-1},f^{k-2}) \ar[r]^{i_*\ \ \ } \ar[d]_{\rho^{-1}} & H_{k-1}\big(f^{k-1},f^{k-1}\setminus \{q\}\big)\\
\mathcal{C}_{k-1}(f) \ar[d]_{\Hat{q}} & H_{k-1}\big(\mathcal{W}_u(q),\mathcal{W}_u(q)^*\big) \ar[u]_{\simeq} \ar[d]^{\simeq} \\
\Z \ar[r]_{Id} & \Z } 
\end{equation} es conmutativo, donde el isomorfismo entre $\Z$ y
$H_{k-1}\big(\mathcal{W}_u(q),\mathcal{W}_u(q)^*\big)$ es asociado a la 
orientaci\'on escogida en $\mathcal{W}_u(q)$, i.e., mapea $\mathcal{O}_q$ en 
$1$. 
\begin{rema}\mbox{}
\begin{enumerate}
\item
Para $p \in \crit _k$, sigue de la definici\'on de $\rho$, que el homeomorfismo
 $f_1$ definido por la composici\'on:
$$
\xymatrix{\Z \ar[rrdd]_{f_1} \ar[r]^{\simeq \ \ \ \ \ \ \ \ \ \ \ \ \ \ \ \ \ } & H_{k}\big(\mathcal{W}_u(p),\mathcal{W}_u(p)^*\big)\ar[r]^{i_*\ \ \ \ \ \ \ \ \ } & 
 \displaystyle \bigoplus_{p_i\in\crit_k} H_{k}\big(\mathcal{W}_u(p_i),\mathcal{W}_u(p_i)^*\big) \ar[d]^{\varphi}\\
& & H_k(f^{k},f^{k-1}) \ar[d]^{\partial _*} \\
& & H_{k-1}(f^{k-1},f^{k-2}) }
$$ mapea $1$ sobre $\partial*\big(\rho(p)\big)$.

\item
Para $q \in \crit _{k-1}$, el homeomorfismo $f_3$ definido por la composici\'on:
\[
\xymatrix{H_{k-1}(f^{k-1},f^{k-2}) \ar[r]^{i_*\ \ \ \ } \ar@/_1.0pc/[rrd]_{f_3} & H_{k-1}\big(f^{k-1},f^{k-1}\setminus \{q\}\big) \ar[r]^{\simeq\ } & H_{k-1}\big(\mathcal{W}_u(q),\mathcal{W}_u(q)^*\big) \ar[d]^{\simeq}\\
& & \Z }
\] coincide con $\Hat{q} \circ \rho^{-1}$. Ver diagrama \eqref{d2}. Luego la 
composici\'on $f_3 \circ f_1$ mapea $1$ sobre
$\hat{q}\big(\rho^{-1}\partial _* \rho (p)\big)$. 
\end{enumerate}
\end{rema}
Dado $a\in \R$ tal que $k-1<a<k$, para cada $p\in \crit_k$, cada 
$q\in \crit_{k-1}$, considere el diagrama:

\begin{equation}\label{d3} \xymatrix{\Z \ar[r]^{Id} \ar[dd]_{f_1} & \Z \ar[d]^{\simeq} \\
& \tilde{H}_{k-1}\left(f^{-1}(a) \cap \mathcal{W}_u(p)\right) \ar[d]^{i_*} \\
H_{k-1}(f^{k-1},f^{k-2}) \ar[dd]_{f_3} & H_{k-1}\left(f^{-1}(a),f^{-1}(a)\setminus \mathcal{W}_s(q)\right) \ar[d]^{\varrho}\\
& H_0\left(f^{-1}(a) \cap \mathcal{W}_s(q)\right) \ar[d]^{\oplus}\\
\Z \ar[r]_{Id} & \Z } 
\end{equation} donde el isomorfismo 
$\Z\simeq\tilde{H}_{k-1}(f^{-1}(a) \cap \mathcal{W}_u(p))$ es asociado a la 
orientaci\'on en la esfera $f^{-1}(a) \cap \mathcal{W}_u(p)$ determinada por 
$\mathcal{O}_p$. i.e., mapea $1$ en $\alpha^{[k-1]}$. La columna izquierda del diagrama \eqref{d3} mapea $1$ sobre
$\hat{q}\big(\rho^{-1}\partial _* \rho (p)\big)$ y la columna de la derecha mapea $1$ en el número de intersección de las esferas $f^{-1}(a) \cap \mathcal{W}_u(p)$ y $f^{-1}(a) \cap \mathcal{W}_s(q)$ en la subvariedad $f^{-1}(a)$, el cual es el enetero $\hat{q}\left(\partial_*(p)\right)$. Por lo tanto, para obtener la conmutatividad del diagrama \eqref{d0}, basta mostrar la conmutatividad del diagrama \eqref{d3}. 

\noindent Para establecer la conmutatividad del diagrama \eqref{d3}, construiremos un diagrama conmutativo equivalente a éste. Para esto fijemos n\'umeros reales $a, b$ tales que $k-1 < a < b < k$ y puntos 
cr\'iticos $p\in \crit_k$, $q\in \crit_{k-1}$. Como conecuencia del Teorema
\ref{moustro}, tenemos un isomorfismo:
\[
H_{k-1}\big(f^{b},f^{b}\setminus \mathcal{W}_s(q)\big)
\xrightarrow{\;\;\varrho\;\;}
H_0\big(f^{b},f^{b}\setminus \mathcal{W}_s(q)\big).
\]adem\'as, una diagrama conmutativo:
\[
\xymatrix{H_{k-1}\big(f^{b},f^b\setminus \mathcal{W}_s(q)\big)  \ar[r]^{\varrho\ \ } & H_0\left(f^{b}\cap \mathcal{W}_s(q)\right)\\
H_{k-1}\left(f^{1}(a),f^{-1}(a)\setminus \mathcal{W}_s(q)\right)
\ar[u]^{\ \ \ i_*}
 \ar[r]_{\ \ \ \ \ \varrho} & H_0\big(f^{-1}(a) \cap \mathcal{W}_s(q)\big) \ar[u]_{i_*} }
\]
Denotamos por $f_2, f_4, f_5$ los homeomorfismos definidos como sigue: 
\[
\xymatrix{%
\Z\ar@{->}@/_2pc/[rr]_{f_2}\ar[r]^-{\simeq}&\tilde{H}_{k-1}\big(f^{-1}(a)\cap
\mathcal{W}_u(p)\big)\ar[r]^-{i_*}&
H_{k-1}\big(f^{b},f^{b}\setminus \mathcal{W}_s(q)\big)} 
\]
\[\xymatrix{%
H_{k-1}\big(f^{b},f^{b}\setminus \mathcal{W}_s(q)\big)\ar@{->}@/_2pc/[rr]_{f_4}\ar[r]^-\varrho
&H_0\big(f^{b}\setminus \mathcal{W}_s(q)\big)
\ar[r]^-{\simeq}&\Z}
\]
\[
\xymatrix{%
H_{k-1}\big(f^{-1}(a),f^{-1}(a)_*\big)\ar@{->}@/_2pc/[rr]_{f_5}
\ar[r]^-{\varrho}&H_0\big(f^{-1}(a) \cap
\mathcal{W}_s(q)\big)\ar[r]^-{\oplus}&\Z}
\] Donde empleamos la notación $f^{-1}(a)_*=f^{-1}(a)\setminus
\mathcal{W}_s(q)$, además el isomorfismo $\Z \simeq
\tilde{H}_{k-1}(f^{-1}(a) \cap \mathcal{W}_u(p))$ es asociado con la 
orientaci\'on en la esfera $f^{-1}(a) \cap \mathcal{W}_u(p)$ determinada por
$\mathcal{O}_p$ y el isomorfismo $\varrho$ tiene signo definido por la 
orientaci\'on transversal en $\mathcal{W}_s(q)$.

Con la notaci\'on y terminolog\'ia anterior, la conmutatividad del 
diagrama \eqref{d3} se establace a partir de la conmutatividad del siguiente 
diagrama: 
\begin{equation}\label{eq:dfinal}
\xymatrix{ &\Z \ar[dd]_{f_2} \ar[dr]^{\simeq \ \ \ \ \ \ } \ar[dl]_{\ \ \ \ \ \ f_1}& \\
H_{k-1}(f^{k-1},f^{k-2})\ar[rd]^{i_*} \ar[rddd]_{f_3}& &\tilde{H}_{k-1}\big(f^{-1}(a)\cap\mathcal{W}_u(p)\big) \ar[d]^{i_*} \ar[dl]_{i_*}\\
&H_{k-1}\big(f^{b},f^{b}\setminus \mathcal{W}_s(q)\big)\ar[dd]_{f_4} & H_{k-1}\big(f^{-1}(a),f^{-1}(a)_*\big) \ar[l]_{i_*\ \ \ } \ar[ddl]^{f_5}\\
& &\\
&\Z& }
\end{equation} Es claro que la conmutatividad del cuadrado en la parte superior derecha es una consecuencia inmediata de la funtorialidad de la homolog\'ia singular. Por lo tanto, resta mostrar la conmutatividad de los triángulos 
que aparecen en el diagrama.
\begin{rema}
\item
Dado $p \in \crit_k$, como $\mathcal W_u(p)$ es un espacio contr\'actil, la 
secuencia larga en homolog\'ia del par 
$\big(\mathcal W_u(p),\mathcal W_u(p)^*\big)$ induce un isomorfismo en 
homolog\'ia:
\[
H_k\big(\mathcal{W}_u(p),\mathcal{W}_u(p)^*\big) \xrightarrow{\;\;\partial _*\;\;}
\tilde{H}_{k-1}\big(\mathcal{W}_u(p)^*\big),
\]  adem\'as  la esfera $\mathcal{W}_u(p)\cap f^{-1}(a)$ es un retrato por deformaci\'on del espacio $\mathcal{W}_u(p)^*$, luego 
la aplicación onclusión induce un isomorfismo en homolog\'ia: 
\[
\widetilde H_{k-1}\big(\mathcal{W}_u(p)\cap f^{-1}(a)\big)
\xrightarrow{\;\;i_*\;\;}\widetilde
H_{k-1}\big(\mathcal{W}_u(p)^*\big).
\] Finalmente, el homeomorfismo en homología
$
H_k\big(f^{b}\big)\xrightarrow{\;\;i_*\;\;}
H_k(f^k\setminus \crit _k)
$ inducido por inclusi\'on es un isomorfismo.
\end{rema}

\begin{lema}\label{a-1}
El siguiente diagrama es conmutativo:
\[
\xymatrix{& \Z \ar[dr]^{f_2} \ar[dl]_{f_1} & \\
H_{k-1}\big(f^{k-1},f^{k-2}\big) \ar[rr]_{i_*} & &
H_{k-1}\big(f^{b},f^{b}\setminus \mathcal{W}_s(q)\big) }
\]
\end{lema}
\proof El diagrama es equivalente con el siguiente diagrama conmutativo:
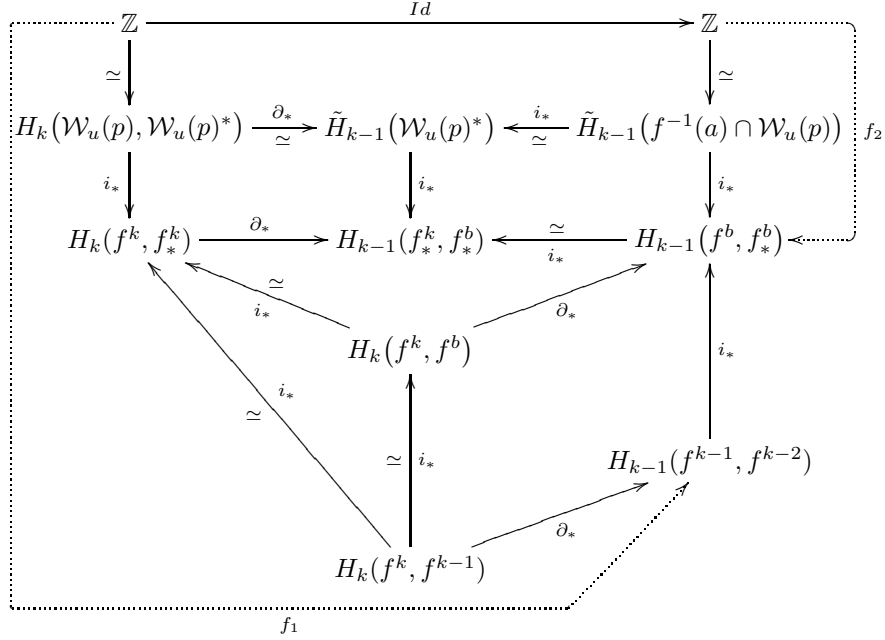
\[
\xymatrix{\Z \ar[rr]^{Id} \ar[d]_{\simeq}
\ar@.`l/1pt[d]-<45pt,0pt>
`[ddddd]-<0pt,15pt>`^u[dddddr]-<-60pt,15pt>_{f_1}[ddddrr]
& & \Z \ar[d]^{\simeq} \ar@.`r/3pt[dd]+<55pt,0pt> `[dd]+<30pt,0pt>^{f_2} [dd]+<30pt,0pt> \\
H_k\big(\mathcal{W}_u(p),\mathcal{W}_u(p)^*\big) \ar[d]_{i_*} \ar[r]^{\ \ \ \ \partial _*}_{\ \ \ \ \simeq} &\tilde{H}_{k-1}\big(\mathcal{W}_u(p)^*\big) \ar[d]^{i_*} &\tilde{H}_{k-1}\big(f^{-1}(a)\cap\mathcal{W}_u(p)\big)\ar[l]_-{\ \ i_*}^-{\ \ \simeq\ }\ar[d]^{i_*}\\
H_k(f^k,f^k_*) \ar[r]^-{\partial _* \ } &H_{k-1}(f^k_*,f^{b}_*) & H_{k-1}\big(f^{b},f^{b}_*\big) \ar[l]_{\simeq\ }^{i_*\ }\\
& H_k\big(f^k,f^{b}\big)\ar[ur]_{\partial _*} \ar[ul]^{i_*}_{\simeq} & \\
&  & H_{k-1}(f^{k-1},f^{k-2}) \ar[uu]_{i_*} \\
& H_k(f^k,f^{k-1}) \ar[uu]_{i_*}^{\simeq} \ar[ur]_{\partial _ *}
\ar[uuul]^{\simeq}_{i_*} &  }
\]
\endproof
donde empleamos la notación $f^{k}_*=f^{k}\setminus \crit_k$, $f^{b}_*=f^{b}\setminus \mathcal{W}_s(p)$.
\begin{lema}\label{a-2}
Con la notaci\'on y la terminolog\'ia anterior. El diagrama:
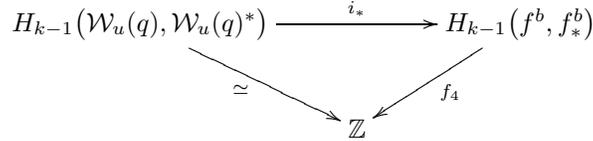
\[
\xymatrix{H_{k-1}\big(\mathcal{W}_u(q),\mathcal{W}_u(q)^*\big)\ar[rr]^-{i_*} \ar[dr]_{\simeq}& & H_{k-1}\big(f^{b},f^{b}_*\big) \ar[dl]^{f_4} \\
& \Z  & }
\] es conmutativo. Adem\'as el isomorfismo
$H_{k-1}\big(\mathcal{W}_u(q),\mathcal{W}_u(q)^*\big) \simeq
\Z$ es asocido a la orientaci\'on transversal en $\mathcal{W}_u(q)$.
\end{lema}
\proof Como $n\_(q) = k-1$,  el conjunto $\mathcal W_u(q)$ siendo homeomorfo a 
$\R^{k-1}$ puede ser identificado con alg\'un conjunto abierto 
de la esfera $S^{k-1}$. Adem\'as, $\mathcal{W}_s(q)\cap
f^{b}$ tiene codimensi\'on $k-1$ en
$f^{b}$. Luego el diagrama:
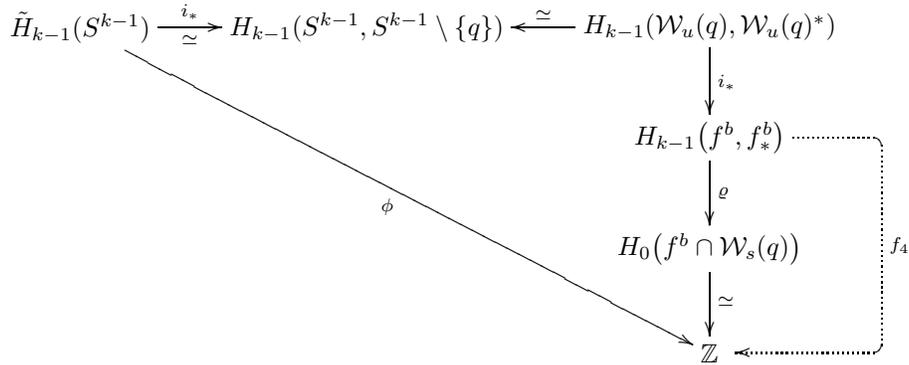
\[
\xymatrix{\tilde{H}_{k-1}(S^{k-1}) \ar[r]_-{\simeq}^-{i_*} \ar[dddrr]_{\phi} & H_{k-1}(S^{k-1},S^{k-1}\setminus \{q\}) & H_{k-1}(\mathcal{W}_u(q),\mathcal{W}_u(q)^*)\ar[l]_-{\simeq} \ar[d]^{i_*}\\
& & H_{k-1}\big(f^{b},f^{b}_*\big) \ar[d]^{\varrho} \ar@.`r/3pt[dd]+<65
pt,0pt> `[dd]+<30pt,0pt>^{f_4} [dd]+<10pt,0pt>\\
& & H_0\big(f^{b}\cap \mathcal{W}_s(q)\big) \ar[d]^{\simeq}\\
& & \Z }
\]
calcula el n\'umero de intersecci\'on de la aplicaci\'on inclusi\'on 
$S^{k-1}\supset\mathcal{W}_u(q)\longrightarrow f^{b}$ con 
$\mathcal{W}_s(q)\cap f^{b}$. Es claro que este n\'umero es $1$,
o sea, que la funci\'on compuesta $f_4 \circ i_*$ lleva $\mathcal{O}_q$ en
$1$.
 \endproof
\begin{lema}\label{a-3}
Con la notaci\'on y la terminolog\'ia anterior. El diagrama:
\[
\xymatrix{H_{k-1}\big(f^{k-1},f^{k-2}\big) \ar[rr]^-{i_*} \ar[dr]_{f_3} & & H_{k-1}\big(f^{b},f^{b}_*\big) \ar[dl]^{f_4}\\
& \Z & }
\] es commutativo.
\end{lema}
\proof El diagrama es equivalente con el siguiente diagrama conmutativo:
\[
\xymatrix{H_{k-1}(f^{k-1},f^{k-2}) \ar[r]^{i_*\ \ } \ar[d]_{i_*} \ar@. `l/3pt[ddd]-<75pt,0pt> `[ddd]-<10pt,0pt>_{f_3} [ddd]-<10pt,0pt> & H_{k-1}\big(f^{b},f^{b}_*\big) \ar[dd]^{\varrho} \ar@.`r/3pt[ddd]+<80pt,0pt> `[ddd]+<10pt,0pt>^{f_4} [ddd]+<10pt,0pt>\\
H_{k-1}(f^{k-1},f^{k-1}\setminus \{q\}) & \\
H_{k-1}\big(\mathcal{W}_u(q),\mathcal{W}_u(q)^*\big) \ar[u]^{i_*}_{\simeq} \ar[uur]_{i_*} \ar[d]_{\simeq} & H_0\big(f^{b}\cap \mathcal{W}_s(q)\big) \ar[d]^{\simeq}\\
\Z \ar[r]_{Id} & \Z }
\]
cuya conmutatividad es garantizada por el lema \ref{a-2} y la funtorialidad de 
la homolog\'ia singular.
\endproof
Finalmente, tenemos 
\begin{lema}\label{a-4}
El siguiente diagrama:
\[
\xymatrix{H_{k-1}\big(f^{b},f^{b}_*\big) \ar[dr]_{f_4} & & H_{k-1}\big(f^{-1}(a),f^{-1}(a)_*\big) \ar[ll]_-{i_* } \ar[dl]^{f_5}\\
& \Z & }
\] es commutativo.
\end{lema}
\proof El diagrama es conmutativo pues es equivalente al siguiente diagrama 
conmutativo:
\[
\xymatrix{H_{k-1}\left(f^{b},f^{b}_*\right) \ar[d]_{\varrho} \ar@.`l/3pt[dd]-<80pt,0pt> `[dd]-<30pt,0pt>_{f_4} [dd]-<10pt,0pt> & H_{k-1}\big(f^{-1}(a),f^{-1}(a)\setminus \mathcal{W}_s(q)\big) \ar[l]_-{\ \ i_*} \ar[d]^{\varrho} \ar@.`r/3pt[dd]+<85pt,0pt> `[dd]+<30pt,0pt>^{f_5} [dd]+<10pt,0pt>\\
H_0\big(f^{b}\cap \mathcal{W}_s(q)\big) \ar[d]^{\simeq} & H_0\big(f^{-1}(a)\cap \mathcal{W}_s(q)\big) \ar[l]_-{\ i_*} \ar[d]^{\oplus}\\
\Z & \Z \ar[l]^{Id}}
\]
\endproof
\end{section}

\end{document}